\documentclass[10pt]{article}
\usepackage{cite}
\usepackage{mathrsfs}
\usepackage{amsfonts}
\usepackage{amsmath}
\usepackage{amsfonts,amssymb,color}
\usepackage{dsfont}
\usepackage{curves}
\usepackage{mathrsfs}
\usepackage{pifont}
\usepackage{amssymb}
\allowdisplaybreaks

\numberwithin{equation}{section}

\date{}

\def\BigRoman{\uppercase\expandafter{\romannumeral\number\count 255 }}
\def\Romannumeral{\afterassignment\BigRoman\count255=}

\begin{document}
\title{Distance spectral radius and perfect matchings in graphs with given fractional property
}
\author{\small Sizhong Zhou\footnote{Corresponding author. E-mail address: zsz\_cumt@163.com (S. Zhou)}\\
\small  School of Science, Jiangsu University of Science and Technology,\\
\small  Zhenjiang, Jiangsu 212100, China\\
}

\maketitle
\begin{abstract}
\noindent A matching in a graph $G$ is a set of independent edges in $G$. A perfect matching in a graph $G$ is a matching which saturates all the vertices of $G$. A fractional perfect matching in a graph
$G$ is a function $h:E(G)\rightarrow [0,1]$ such that $\sum\limits_{e\in E_G(v)}h(e)=1$ for every $v\in V(G)$, where $E_G(v)$ is the set of edges incident to $v$ in $G$. Clearly, the existence of a
fractional perfect matching in a graph is a necessary condition for the graph to possess a perfect matching. Let $G$ be a $k$-connected graph of even order $n$ with a fractional perfect matching, where $k$
is a positive integer. We denote by $\mu(G)$ the distance spectral radius of $G$. In this paper, we prove that if $n\geq8k+6$ and $\mu(G)\leq\mu(K_k\vee(kK_1\cup K_3\cup K_{n-2k-3}))$, then $G$ contains a
perfect matching unless $G=K_k\vee(kK_1\cup K_3\cup K_{n-2k-3})$.
\\
\begin{flushleft}

(2020) Mathematics Subject Classification: 05C50, 05C70
\end{flushleft}
\end{abstract}

\section{Introduction}

Let $G$ be a simple graph with vertex set $V(G)$ and edge set $E(G)$. The order of $G$ is denoted by $|V(G)|=n$. For a subset $S\subseteq V(G)$, $G[S]$ and $|S|$ denote the subgraph of $G$ induced by $S$
and the size of $S$, respectively. We write $G-S$ for $G[V(G)\setminus S]$. The number of odd components and the number of isolated vertices in $G$ are denoted by $o(G)$ and $i(G)$, respectively. We denote
by $K_n$ the complete graph of order $n$. Let $G_1$ and $G_2$ be two vertex-disjoint graphs. We use $G_1\cup G_2$ to denote the union of $G_1$ and $G_2$. We write $G_1\vee G_2$ for the new graph, called
the join of $G_1$ and $G_2$, constructed from $G_1\cup G_2$ by joining every vertex of $G_1$ to all the vertices of $G_2$. Let $tG$ denote the union of $t$ copies of $G$.

Given a connected graph $G$ with $V(G)=\{v_1,v_2,\ldots,v_n\}$, the distance between $v_i$ and $v_j$ in $G$, denoted by $d_{ij}$, is defined as the length of the shortest path from $v_i$ to $v_j$ in $G$.
The adjacency matrix of $G$, denoted by $A(G)$, is the matrix $A(G)=(a_{ij})_{n\times n}$, where $a_{ij}=1$ if $v_iv_j$ is an edge of $G$, and 0 otherwise. The adjacency spectral radius of $G$, denoted by
$\rho(G)$, is the largest eigenvalue of its adjacency matrix $A(G)$. The distance matrix of $G$, denoted by $\mathcal{D}(G)$, is defined to be the matrix $\mathcal{D}(G)=(d_{ij})_{n\times n}$. The distance
spectral radius of $G$, denoted by $\mu(G)$, is defined as the largest eigenvalue of the distance matrix $\mathcal{D}(G)$ of $G$. Some properties of the spectral radius of $G$ can be referred to
\cite{Em,Od,Ws1,Ws2,ZW,Za,ZZZL,Zs2,Zs3,FL}.

A matching in $G$ is a set of independent edges in $G$. A perfect matching in $G$ is a matching which saturates all the vertices of $G$. A fractional perfect matching in $G$ is a function
$h:E(G)\rightarrow [0,1]$ such that $\sum\limits_{e\in E_G(v)}h(e)=1$ for every $v\in V(G)$, where $E_G(v)$ is the set of edges incident to $v$ in $G$.

The relationships between graph factors and the spectral radii of various graph matrices have been studied by many scholars. O \cite{Os} we established a lower bound for the adjacency spectral radius in a
connected graph of order $n$ to guarantee the existence of a perfect matching. Zhao, Huang and Wang \cite{ZHW} provided a sufficient condition for the existence of a perfect matching in a connected graph
based on the generalized adjacency spectral radius. Jia, Fan and Liu \cite{JFL} presented an adjacency spectral radius condition for a connected graph with a fractional perfect matching to has a perfect
matching. Zhou, Bian and Wu \cite{ZBW} showed a sufficient conditions for a connected graph with the minimum degree $\delta$ to contain an even factor based its adjacency spectral radius. Wu \cite{Wc}
posed an adjacency spectral radius condition for the existence of a spanning tree with leaf degree at most $k$ in a connected graph. Zhou, Bian and Sun \cite{ZBS}, Pan and Zhou \cite{PZ} gave some sufficient
conditions in terms of size, generalized adjacency spectral radius, distance spectral radius and distance signless Laplacian spectral radius to guarantee the existence of path-factors in isolated tough graphs.
Zhou \cite{Zs1} arose two sufficient conditions for a connected bipartite graph to possess a $[1,k]$-factor with respect to the adjacency spectral radius and the distance spectral radius. Miao and Li \cite{ML}
showed sufficient conditions based on the adjacency spectral radius and the distance spectral radius to guarantee the existence of a star-factor in a connected graph.

Zhang and Lin \cite{ZL} provided a distance spectral radius condition for a connected graph with a perfect matching.

\medskip

\noindent{\textbf{Theorem 1.1}} (Zhang and Lin \cite{ZL}). Let $G$ be a connected graph of even order $n\geq4$.

(\romannumeral1) For $n\leq8$, if $\mu(G)\leq\mu(K_{\frac{n}{2}-1}\vee(\frac{n}{2}+1)K_1)$, then $G$ has a perfect matching unless $G=K_{\frac{n}{2}-1}\vee(\frac{n}{2}+1)K_1$.

(\romannumeral2) For $n\geq10$, if $\mu(G)\leq\mu(K_1\vee(K_{n-3}\cup2K_1))$, then $G$ has a perfect matching unless $G=K_1\vee(K_{n-3}\cup2K_1)$.

\medskip

Notice that the existence of a fractional perfect matching in a graph is a necessary condition for the graph to possess a perfect matching. Together with Theorem 1.1, it is natural that we put forward the
following problem.

\medskip

\noindent{\textbf{Problem 1.2.}} What is a tight distance spectral radius condition which guarantees a graph with a fractional perfect matching to possess a perfect matching?

\medskip

Focusing on Problem 1.2, we verify the following result.

\medskip

\noindent{\textbf{Theorem 1.3.}} Let $k$ be a positive integer and let $G$ be a $k$-connected graph of even order $n\geq8k+6$ with a fractional perfect matching. If
$$
\mu(G)\leq\mu(K_k\vee(kK_1\cup K_3\cup K_{n-2k-3})),
$$
then $G$ has a perfect matching unless $G=K_k\vee(kK_1\cup K_3\cup K_{n-2k-3})$.

\medskip

If $k=1$ in Theorem 1.3, then we possess the following corollary.

\medskip

\noindent{\textbf{Corollary 1.4.}} Let $G$ be a connected graph of even order $n\geq14$ with a fractional perfect matching. If
$$
\mu(G)\leq\mu(K_1\vee(K_1\cup K_3\cup K_{n-5})),
$$
then $G$ has a perfect matching unless $G=K_1\vee(K_1\cup K_3\cup K_{n-5})$.

\medskip

One easily checks that $\mu(K_1\vee(K_1\cup K_3\cup K_{n-5}))>\mu(K_1\vee(K_{n-3}\cup2K_1))$ for $n\geq14$. Obviously, the distance spectral radius condition in Corollary 1.4 is better than that of Theorem 1.1
when $n\geq14$. Hence, Theorem 1.3 is an improvement and generalization of Theorem 1.1.

\section{Preliminary lemmas}

In this section, we introduce several necessary lemmas to verify our main result.

\medskip

\noindent{\textbf{Lemma 2.1}} (Tutte \cite{Tt}). A graph $G$ contains a perfect matching if and only if
$$
o(G-S)\leq|S|
$$
for every vertex subset $S$ of $G$.

\medskip

\noindent{\textbf{Lemma 2.2}} (Scheinerman and Ullman \cite{SU}). A graph $G$ contains a fractional perfect matching if and only if
$$
i(G-S)\leq|S|
$$
for every vertex subset $S$ of $G$.

\medskip

\noindent{\textbf{Lemma 2.3}} (Min\'c \cite{Mn}). Let $G$ be a connected graph with $u,v\in V(G)$ and $uv\notin E(G)$. Then
$$
\mu(G+uv)<\mu(G).
$$

\medskip

The Wiener index of a connected graph $G$ of order $n$ is defined by $W(G)=\sum\limits_{i<j}d_{ij}$. Based on the Rayleigh quotient \cite{HJ}, the following result is easily obtained.

\medskip

\noindent{\textbf{Lemma 2.4.}} Let $G$ be a connected graph of order $n$. Then
$$
\mu(G)=\max_{X\neq\mathbf{0}}\frac{X^{T}\mathcal{D}(G)X}{X^{T}X}\geq\frac{\mathbf{1}^{T}\mathcal{D}(G)\mathbf{1}}{\mathbf{1}^{T}\mathbf{1}}=\frac{2W(G)}{n},
$$
where $\mathbf{1}=(1,1,\ldots,1)^{T}$.

\medskip

\noindent{\textbf{Lemma 2.5}} (Hu and Zhang \cite{HZ}). Let $q$ and $s$ be positive integers with $q\geq s+2$. Let $n_1,\ldots,n_q$ be odd integers with $n_q\geq\cdots\geq n_1\geq1$, $n_{s+1}\geq3$,
$\sum\limits_{i=1}^{q}n_i=n-s$ and $n_q<n-3q+s+3$. Then
$$
\mu(K_s\vee(K_{n_1}\cup\cdots\cup K_{n_q}))>\mu(K_s\vee(sK_1\cup(q-s-1)K_3\cup K_{n-3q+s+3})).
$$

\medskip

Let $M$ denote a real $n\times n$ matrix whose rows and columns are indexed by the set $\mathcal{N}=\{1,2,\ldots,n\}$. Given a partition
$$
\pi:\mathcal{N}=\mathcal{N}_1\cup\mathcal{N}_2\cup\cdots\cup\mathcal{N}_t
$$
of the index set $\mathcal{N}$. Based on the partition $\pi$, the matrix $M$ can be written in block form as
\begin{align*}
M=\left(
  \begin{array}{cccc}
    M_{11} & M_{12} & \cdots & M_{1t}\\
    M_{21} & M_{22} & \cdots & M_{2t}\\
    \vdots & \vdots & \ddots & \vdots\\
    M_{t1} & M_{t2} & \cdots & M_{tt}\\
  \end{array}
\right),
\end{align*}
where the block $M_{ij}$ denotes the $n_i\times n_j$ matrix for $1\leq i,j\leq t$. Let $m_{ij}$ denote the average row sum of the block $M_{ij}$ of $M$ for $1\leq i,j\leq t$. Then the matrix
$M_{\pi}=(m_{ij})_{t\times t}$ is called the quotient matrix of $M$ corresponding to the partition $\pi$. The partition $\pi$ is called equitable if every block $M_{ij}$ of $M$ has constant row sum for
$1\leq i,j\leq t$.

\medskip

\noindent{\textbf{Lemma 2.6}} (Brouwer and Haemers \cite{BH}, You, Yang, So and Xi \cite{YYSX}). Let $M$ be a real $n\times n$ matrix with an equitable partition $\pi$, and let $M_{\pi}$ be the corresponding
quotient matrix. Then the eigenvalues of $M_{\pi}$ are eigenvalues of $M$. Furthermore, if $M$ is nonnegative and irreducible, then the largest eigenvalues of $M$ and $M_{\pi}$ are equal.

\medskip

\section{The proof of Theorem 1.3}

\noindent{\it Proof of Theorem 1.3.} Suppose, to the contrary, that a $k$-connected graph $G$ with a fractional perfect matching contains no perfect matching. By virtue of Lemma 2.1, we conclude
$o(G-S)\geq|S|+1$ for some subset $S$ of $V(G)$. Since $n$ is even, $o(G-S)$ and $|S|$ have the same parity. Thus, we possess $o(G-S)\geq|S|+2$. We are to prove $|S|\geq k$. In fact, if $|S|\leq k-1$,
then $1\geq o(G-S)\geq|S|+2\geq2$ due to $G$ being $k$-connected. This is impossible. Hence, $|S|\geq k$.

Let $|S|=s$ and $o(G-S)=q$. Then we have $q\geq s+2$. The $q$ odd components in $G-S$ are denoted by $O_1,O_2,\ldots,O_q$, where $|O_i|=n_i$ for $1\leq i\leq q$. Without loss of generality, we let
$n_q\geq n_{q-1}\geq\cdots\geq n_2\geq n_1$.

\medskip

\noindent{\bf Claim 1.} $n_{s+1}\geq3$.

\noindent{\it Proof.} Assume that $n_{s+1}=1$. Then we easily see that $n_1=n_2=\cdots=n_{s+1}=1$ and $i(G-S)\geq s+1$. Notice that $G$ has a fractional perfect matching. In terms of Lemma 2.2, we
deduce $i(G-S)\leq s$ for each vertex subset $S$ of $G$ with $s\geq k$. From the above discussion, we obtain $s+1\leq i(G-S)\leq s$, a contradiction. Hence, we infer $n_{s+1}\geq3$. This completes
the proof of Claim 1. \hfill $\Box$

\medskip

In terms of Claim 1, we obtain $n_i\geq3$ for $s+1\leq i\leq q$ and $n_j\geq1$ for $1\leq j\leq s$. Obviously, $G$ is a spanning subgraph of $G_1=K_s\vee(K_{n_1}\cup\cdots\cup K_{n_s}\cup K_{n_{s+1}}\cup K_{n'_{s+2}})$
for some positive odd integers $n'_{s+2}\geq n_{s+1}\geq n_s\geq\cdots\geq n_1\geq1$ with $n_{s+1}\geq3$ and $\sum\limits_{i=1}^{s+1}n_i+n'_{s+2}=n-s$. According to Lemma 2.3, we have
\begin{align}\label{eq:3.1}
\mu(G)\geq\mu(G_1),
\end{align}
where the equality follows if and only if $G=G_1$.

Let $G_2=K_s\vee(sK_1\cup K_3\cup K_{n-2s-3})$, where $n\geq2s+6$. Based on Lemma 2.5, we get
\begin{align}\label{eq:3.2}
\mu(G_1)\geq\mu(G_2),
\end{align}
where the equality occurs if and only if $(n_1,\ldots,n_{s},n_{s+1},n'_{s+2})=(1,\ldots,1,3,n-2s-3)$. In what follows, we proceed to prove this theorem in terms of two possible cases.

\noindent{\bf Case 1.} $s=k$.

In this case, $G_2=K_k\vee(kK_1\cup K_3\cup K_{n-2k-3})$. By virtue of \eqref{eq:3.1} and \eqref{eq:3.2}, we conclude
$$
\mu(G)\geq\mu(K_k\vee(kK_1\cup K_3\cup K_{n-2k-3})),
$$
with equality holding if and only if $G=K_k\vee(kK_1\cup K_3\cup K_{n-2k-3})$. Observe that $K_k\vee(kK_1\cup K_3\cup K_{n-2k-3})$ contains no perfect matching. Thus, we get a contradiction.

\noindent{\bf Case 2.} $s\geq k+1$.

Recall that $G_2=K_s\vee(sK_1\cup K_3\cup K_{n-2s-3})$. The quotient matrix of $\mathcal{D}(G_2)$ based on the partition $V(G_2)=V(K_s)\cup V(sK_1)\cup V(K_3)\cup V(K_{n-2s-3})$ is equal to
\begin{align*}
B_2=\left(
  \begin{array}{cccc}
  s-1 & s & 3 & n-2s-3\\
  s & 2s-2 & 6 & 2n-4s-6\\
  s & 2s & 2 & 2n-4s-6\\
  s & 2s & 6 & n-2s-4\\
  \end{array}
\right).
\end{align*}
Based on a direct computation, the characteristic polynomial of $B_2$ is
\begin{align*}
f_{B_2}(x)=&x^{4}-(n+s-5)x^{3}-((2s+13)n-5s^{2}-16s-36)x^{2}\\
&+((s^{2}-7s-32)n-2s^{3}+16s^{2}+62s+88)x\\
&+(4s^{2}-2s-20)n-8s^{3}-4s^{2}+36s+56.
\end{align*}
In view of Lemma 2.6 and the equitable partition $V(G_2)=V(K_s)\cup V(sK_1)\cup V(K_3)\cup V(K_{n-2s-3})$, we know that the largest root of $f_{B_2}(x)=0$ is equal to $\mu(G_2)$.

Let $G_*=K_k\vee(kK_1\cup K_3\cup K_{n-2k-3})$. The quotient matrix of $\mathcal{D}(G_*)$ corresponding to the partition $V(G_*)=V(K_k)\cup V(kK_1)\cup V(K_3)\cup V(K_{n-2k-3})$ equals
\begin{align*}
B_*=\left(
  \begin{array}{cccc}
  k-1 & k & 3 & n-2k-3\\
  k & 2k-2 & 6 & 2n-4k-6\\
  k & 2k & 2 & 2n-4k-6\\
  k & 2k & 6 & n-2k-4\\
  \end{array}
\right),
\end{align*}
and the characteristic polynomial of $B_*$ is
\begin{align*}
f_{B_*}(x)=&x^{4}-(n+k-5)x^{3}-((2k+13)n-5k^{2}-16k-36)x^{2}\\
&+((k^{2}-7k-32)n-2k^{3}+16k^{2}+62k+88)x\\
&+(4k^{2}-2k-20)n-8k^{3}-4k^{2}+36k+56.
\end{align*}
Based on Lemma 2.6 and the equitable partition $V(G_*)=V(K_k)\cup V(kK_1)\cup V(K_3)\cup V(K_{n-2k-3})$, we see that the largest root of $f_{B_*}(x)=0$ equals $\mu(G_*)$.

Write $\mu=\mu(G_*)$. Recall that $G_*=K_k\vee(kK_1\cup K_3\cup K_{n-2k-3})$. According to Lemma 2.4 and $n\geq8k+6$, we obtain
\begin{align}\label{eq:3.3}
\mu=&\mu(K_k\vee(kK_1\cup K_3\cup K_{n-2k-3}))\nonumber\\
\geq&\frac{2W(K_k\vee(kK_1\cup K_3\cup K_{n-2k-3}))}{n}\nonumber\\
=&\frac{n^{2}+(2k+5)n-3k^{2}-13k-18}{n}\nonumber\\
\geq&\frac{n^{2}+(k+3)n+(k+2)(8k+6)-3k^{2}-13k-18}{n}\nonumber\\
=&\frac{n^{2}+(k+3)n+5k^{2}+9k-6}{n}\nonumber\\
>&n+k+3.
\end{align}

Notice that $f_{B_*}(\mu)=0$. By plugging the value $\mu$ into $x$ of $f_{B_2}(x)-f_{B_*}(x)$, we obtain
\begin{align}\label{eq:3.4}
f_{B_2}(\mu)=f_{B_2}(\mu)-f_{B_*}(\mu)=(s-k)\varphi(s),
\end{align}
where $\varphi(s)=-(2\mu+8)s^{2}+(5\mu^{2}+(n-2k+16)\mu+4n-8k-4)s-\mu^{3}-(2n-5k-16)\mu^{2}+((k-7)n-2k^{2}+16k+62)\mu+(4k-2)n-8k^{2}-4k+36$. Note that
$$
\frac{5\mu^{2}+(n-2k+16)\mu+4n-8k-4}{2(2\mu+8)}>\frac{n-6}{2}\geq s
$$
due to \eqref{eq:3.3} and $n\geq2s+6$. Then $\varphi(s)$ is strictly increasing for $s\leq\frac{n-6}{2}$ and
\begin{align}\label{eq:3.5}
\varphi(s)\leq&\varphi\Big(\frac{n-6}{2}\Big)\nonumber\\
=&-(2\mu+8)\Big(\frac{n-6}{2}\Big)^{2}+(5\mu^{2}+(n-2k+16)\mu+4n-8k-4)\Big(\frac{n-6}{2}\Big)\nonumber\\
&-\mu^{3}-(2n-5k-16)\mu^{2}+((k-7)n-2k^{2}+16k+62)\mu\nonumber\\
&+(4k-2)n-8k^{2}-4k+36\nonumber\\
=&\frac{1}{2}\Big(-2\mu^{3}+(n+10k+2)\mu^{2}+(8n-4k^{2}+44k-8)\mu+16n-16k^{2}+40k-48\Big).
\end{align}
Let $g(\mu)=-2\mu^{3}+(n+10k+2)\mu^{2}+(8n-4k^{2}+44k-8)\mu+16n-16k^{2}+40k-48$. Note that $\mu>n+k+3$ due to \eqref{eq:3.3}. Next, we shall prove that $g(\mu)<0$ for $\mu\geq n+k+3$.

In fact, $g'(\mu)=-6\mu^{2}+2(n+10k+2)\mu+8n-4k^{2}+44k-8$, and the symmetry axis of $g'(\mu)$ is $\mu=\frac{n+10k+2}{6}<n+k+3$. For $\mu\geq n+k+3$, we possess
\begin{align*}
g'(\mu)\leq&g'(n+k+3)\\
=&-6(n+k+3)^{2}+2(n+10k+2)(n+k+3)+8n-4k^{2}+44k-8\\
=&-4n^{2}+(10k-18)n+10k^{2}+72k-50\\
\leq&-4(8k+6)^{2}+(10k-18)(8k+6)+10k^{2}+72k-50 \ \ \ \ \ (\mbox{since} \ n\geq8k+6)\\
=&-166k^{2}-396k-302\\
<&0 \ \ \ \ \ (\mbox{since} \ k\geq1).
\end{align*}
This implies that $g(\mu)$ is strictly decreasing in the interval $[n+k+3,+\infty)$. For $\mu\geq n+k+3$, we deduce
\begin{align}\label{eq:3.6}
g(\mu)\leq&g(n+k+3)\nonumber\\
=&-2(n+k+3)^{3}+(n+10k+2)(n+k+3)^{2}\nonumber\\
&+(8n-4k^{2}+44k-8)(n+k+3)+16n-16k^{2}+40k-48\nonumber\\
=&-n^{3}+(6k-2)n^{2}+(11k^{2}+86k-1)n+4k^{3}+60k^{2}+212k-108.
\end{align}
Let $h(n)=-n^{3}+(6k-2)n^{2}+(11k^{2}+86k-1)n+4k^{3}+60k^{2}+212k-108$. It follows from $k\geq1$ and $n\geq8k+6$ that
\begin{align*}
h'(\mu)=&-3n^{2}+2(6k-2)n+11k^{2}+86k-1\\
\leq&-3(8k+6)^{2}+2(6k-2)(8k+6)+11k^{2}+86k-1\\
=&-85k^{2}-162k-133\\
<&0,
\end{align*}
which implies that $h(n)$ is strictly decreasing with respect to $n\geq8k+6$. Hence, we obtain
\begin{align}\label{eq:3.7}
h(n)\leq&h(8k+6)\nonumber\\
=&-(8k+6)^{3}+(6k-2)(8k+6)^{2}+(11k^{2}+86k-1)(8k+6)\nonumber\\
&+4k^{3}+60k^{2}+212k-108\nonumber\\
=&-36k^{3}+110k^{2}-120k-402\nonumber\\
<&0 \ \ \ \ \ (\mbox{since} \ k\geq1).
\end{align}
From \eqref{eq:3.6} and \eqref{eq:3.7}, we infer $g(\mu)<0$ for $\mu\geq n+k+3$. Combining this with \eqref{eq:3.4}, \eqref{eq:3.5} and $s\geq k+1$, we deduce
$$
f_{B_2}(\mu)=(s-k)\varphi(s)\leq\frac{1}{2}(s-k)g(\mu)<0,
$$
which implies $\mu(G_2)>\mu=\mu(G_*)$. Together with \eqref{eq:3.1} and \eqref{eq:3.2}, we conclude
$$
\mu(G)\geq\mu(G_1)\geq\mu(G_2)>\mu(G_*)=\mu(K_k\vee(kK_1\cup K_3\cup K_{n-2k-3})),
$$
which contradicts $\mu(G)\leq\mu(K_k\vee(kK_1\cup K_3\cup K_{n-2k-3}))$. Theorem 1.3 is proved. \hfill $\Box$

\section*{Declaration of competing interest}

The author declares that he has no known competing financial interests or personal relationships that could have appeared to influence the work reported in this paper.

\section*{Data availability}

No data was used for the research described in the article.

\section*{Acknowledgments}

This work was supported by the Natural Science Foundation of Jiangsu Province (Grant No. BK20241949). Project ZR2023MA078 supported by Shandong Provincial Natural Science Foundation.

\end{document}